\documentclass[12pt,a4paper]{amsart}
\usepackage{amssymb,amscd}
\newtheorem{theorem}{Theorem}
\newtheorem{corollary}{Corollary}
\newtheorem{proposition}{Proposition}
\newtheorem{lemma}{Lemma}

\newcommand\bP{{\mathbb P}}
\newcommand\bQ{{\mathbb Q}}
\newcommand\bZ{{\mathbb Z}}

\newcommand\cF{{\mathcal F}}
\newcommand\cG{{\mathcal G}}
\newcommand\cL{{\mathcal L}}
\newcommand\cN{{\mathcal N}}
\newcommand\cO{{\mathcal O}}
\newcommand\cT{{\mathcal T}}
\newcommand\cV{{\mathcal V}}

\newcommand\fb{{\mathfrak b}}
\newcommand\fk{{\mathfrak k}}
\newcommand\fg{{\mathfrak g}}
\newcommand\ft{{\mathfrak t}}
\newcommand\fu{{\mathfrak u}}

\newcommand\diag{\operatorname{diag}}
\newcommand\id{\operatorname{id}}
\newcommand\rk{\operatorname{rk}}

\newcommand\Ind{\operatorname{Ind}}

\title{Construction of equivariant vector bundles} 

\author{Michel~Brion}
\address{Universit\'e de Grenoble I\\
D\'epartement de Math\'ematiques\\
Institut Fourier, UMR 5582 du CNRS\\
38402 Saint-Martin d'H\`eres Cedex, France}
\email{Michel.Brion@ujf-grenoble.fr}

\begin{document}

\begin{abstract}
Let $X$ be the wonderful compactification of a complex adjoint
symmetric space $G/K$ such that $\rk(G/K)=\rk(G)-\rk(K)$. We show 
how to extend equivariant vector bundles on $G/K$ to equivariant
vector bundles on $X$, generated by their global sections and having
trivial higher cohomology groups. This relies on a geometric
construction of equivariant vector bundles in the setting of
varieties with reductive group action and ``multiplicity-free''
subvarieties. 
\end{abstract}

\maketitle

\section*{Introduction}

Consider a homogeneous space $X_0 = G/K$, where $G$ is a complex 
linear algebraic group and $K$ is a complex algebraic subgroup.
Let $\cV_0$ be a $G$-equivariant vector bundle on $X_0$
and let $X$ be a $G$-equivariant compactification of $G/K$. 
Does $\cV_0$ extend to an equivariant vector bundle on $X$ ? 

\medskip

This question arises naturally when studying equivariant vector 
bundles on quasi-homogeneous varieties. It was raised by Kostant 
in the setting where $X_0$ is an adjoint symmetric space and $X$
is its wonderful compactification (as introduced in \cite{DePr83}), 
in view of applications to representation theory of real reductive
groups. In fact, Kostant asked for a canonical extension; see Syu
Kato's paper \cite{Ka02} for a precise formulation of his question,
and a positive answer in the case of an adjoint semisimple
group $K$ regarded as the symmetric space $K\times K/\diag(K)$.

\medskip

Returning to the general setting, recall that the equivariant vector
bundles on $G/K$ are in bijection with the $K$-modules $M$, via 
$M \mapsto \cL_{G/K}(M)$. Further, if the subgroup $K$ is reductive, 
then any $K$-module is semisimple, and any simple $K$-module $M$ is 
a quotient of some simple $G$-module $V$. Thus, $\cL_{G/K}(M)$ 
is generated by its subspace $V$ of global sections. It follows 
easily that the sheaf of sections of $\cL_{G/K}(M)$ extends to 
a unique $G$-linearized torsion-free, coherent sheaf on $X$, 
generated by its subspace $V$ of global sections; see Lemma 
\ref{unique}. So we may ask : Is this extension is locally free 
for a suitable choice of $V$ ?

\medskip

The answer to this refined question may be negative, as 
$\cL_{G/K}(M)$ may admit no extension as an equivariant vector bundle. 
For example, if $G/K = GL_n/O_n$ is the symmetric 
space of nondegenerate quadratic forms in $n$ variables, and $X$ 
is the space of all quadratic forms, then only the trivial 
equivariant vector bundles (i.e., those associated with $G$-modules)
extend to equivariant vector bundles on $X$. Indeed, any equivariant 
vector bundle on $X$ is trivial by \cite[Prop.10.1]{BaHa85}. But 
clearly, our refined question has a positive answer for a certain 
equivariant compactification $X(M,V)$, and hence for all the 
compactifications which dominate $X(M,V)$.

\medskip

In this article, we obtain an affirmative answer for the wonderful
compactifications of a class of adjoint symmetric spaces. The rank 
of any such space satisfies the inequality 
$\rk(G/K) \ge \rk(G)-\rk(K)$; we say that $G/K$ is of minimal rank 
if the equality holds. We may now formulate a slightly weaker form 
of our main result. 

\begin{theorem}\label{main}
Let $X$ be the wonderful compactification of a complex adjoint
symmetric space $G/K$ of minimal rank, and let $M$ be a simple 
$K$-module. 

Then $M$ is the quotient of some simple $G$-module $V$ such that 
the highest weight line of $V$ projects onto the highest weight 
line of $M$. 

For any such $V$, the bundle $\cL_{G/K}(M)$ extends to a unique 
$G$-equivariant vector bundle on $X$, generated by its space 
$V$ of global sections.
\end{theorem}

Here is an outline of the proof. By the Borel-Weil theorem, 
$M$ is the space of global sections of some $K$-equivariant line 
bundle $\cL_{K/B_K}(\mu)$, where $B_K$ is a Borel subgroup of $K$, 
and $\mu$ is a character of $B_K$ identified with the 
corresponding one-dimensional $B_K$-module. In other words, 
$\cL_{G/K}(M)$ is the direct image of the equivariant 
line bundle $\cL_{G/B_K}(\mu)$ under the equivariant fibration 
$$
\pi_0:G/B_K \to G/K
$$
with fiber being the flag variety $K/B_K$. To obtain the desired
extension, we will construct an equivariant compactification of
$\pi_0$ over $X$, where $\cL_{G/B_K}(\mu)$ extends to an equivariant
line bundle.

\medskip

We may find a Borel subgroup $B$ of $G$ such that $B\cap K = B_K$. 
Then the orbit $Y_0 := B/B_K$ is closed in $G/K$, since $K/B_K$ is
closed in $G/B$. Let $Y$ be the closure of $Y_0$ in $X$, this is a
$B$-stable subvariety. Thus, we may form the ``induced'' $G$-variety
$G\times^B Y$, a compactification of $G\times^B B/B_K \cong G/B_K$.
This variety is provided with $G$-equivariant morphisms 
$$
\pi: G\times^B Y \to X
$$ 
(a compactification of $\pi_0$) and
$$
f: G\times^B Y \to G/B
$$ 
(a compactification of the natural map $f_0: G/B_K \to G/B$).
Further, for any character $\lambda$ of $B$, we have an equivariant
line bundle $f^*\cL_{G/B}(\lambda)$ on $G\times^B Y$; it extends
$\cL_{G/B_K}(\mu)$ if and only if $\lambda$ extends $\mu$. Under 
this assumption, 
$$
\Phi(\lambda) := \pi_*(f^*\cL_{G/B}(\lambda))
$$
is a torsion-free, coherent $G$-linearized sheaf on $X$ which 
extends $\cL_{G/K}(M)$. If, in addition, $\lambda$ is dominant,
then the space $V$ of global sections of $\Phi(\lambda)$ is 
a simple $G$-module, and the restriction to the fiber $K/B_K$ yields 
a surjective map $V\to M$ sending the highest weight line of $V$
onto that of $M$.

\medskip

To check that $\Phi(\lambda)$ is locally free and generated by $V$,
we will apply the theorem on cohomology and base change to the 
morphism $\pi$. Indeed, under the assumptions of Theorem \ref{main}, 
$\pi$ turns out to be flat with reduced fibers. In fact, the fibers 
of $\pi$ realise a flat degeneration of the flag variety $K/B_K$ 
to a union of Schubert varieties in the larger flag variety $G/B$. 
Further, the dominant weight $\mu$ of $K$ turns out to extend to 
a dominant weight $\lambda$ of $G$. So we may complete the proof 
by using known properties of the restriction of $\cL_{G/B}(\lambda)$ 
to unions of Schubert varieties, together with the semicontinuity 
theorem for $\pi$.

\medskip

Several steps of this argument make sense in much greater generality,
and this is how they are presented here. In Section 1, we consider a 
$G$-variety $X$ and a $B$-stable subvariety $Y$, and we study the
geometry of the corresponding morphism $\pi:G \times^B Y \to X$. Our
main result asserts that $\pi$ is flat with reduced fibers whenever
$X$ is complete and nonsingular, and $Y$ is multiplicity-free in the
sense of \cite{Br03}.

\medskip

In Section 2, we study the functors from $B$-modules to $G$-linearized
sheaves on $X$ defined by $M\mapsto R^i\pi_*(f^*\cL_{G/B}(M))$. 
From the results of Section 1, we deduce several relations of these 
functors to the Joseph functors of $B$-module theory (\cite{Po89}, 
\cite{vdK93}). Also, we show that the sheaf 
$\Phi(\lambda) := \pi_*(f^*\cL_{G/B}(\lambda))$ 
is locally free and generated by its global sections, if $\lambda$ 
is a dominant weight, $X$ is complete and nonsingular, and 
$Y$ is multiplicity-free. 

\medskip

The final Section 3 is devoted to symmetric spaces of minimal rank. 
We show that any such symmetric space $G/K$ contains a unique closed
$B$-orbit; further, the closure $Y$ of this orbit in the wonderful 
compactification $X$ is multiplicity-free, and any dominant weight 
of $K$ extends to a dominant weight of $G$. Also, the subspace of 
the rational Grothendieck group $K(X)_{\bQ}$ generated by the 
classes $[\Phi(\lambda)]$, where $\lambda$ is dominant, admits a basis
consisting of the classes of the structure sheaves of closures in $X$ 
of $B$-orbits in $G/K$. As a consequence, the $[\Phi(\lambda)]$
do not generate the whole Grothendieck group; this is discussed in 
more detail in 3.5.

\medskip

Among the adjoint symmetric spaces, those of minimal rank form a rather
restricted class: it consists of the products of homogeneous spaces
$K\times K/\diag(K)$ (where the group $K$ is simple),
$PSL_{2n}/PSp_{2n}$, $PSO_{2n}/PSO_{2n-1}$, and $E_6/F_4$. The
argument of Theorem \ref{main} extends to further examples of 
symmetric spaces, e.g., to $GL_{m+n}/GL_m\times GL_n$ and
$SO_{2n}/GL_n$; indeed, these spaces $G/K$ contain several closed
$B$-orbits, but all of them are multiplicity-free, and any dominant
weight of $K$ extends to a dominant weight of $G$. This argument also
extends to the spherical homogeneous spaces of minimal rank, as
introduced and classified by Ressayre in \cite{Re04}. 

\medskip

It would be interesting to describe the categories of equivariant 
vector bundles on wonderful compactifications of symmetric spaces 
of minimal rank, generalizing Syu Kato's description \cite{Ka02} 
in the case of group compactifications, and to obtain an intrinsic 
characterization of those bundles constructed here. A related open
question is to find a presentation of the rational cohomology rings 
of such compactifications. These rings are generated (as rational
vector spaces) by Chern classes of equivariant vector bundles, but
no minimal set of multiplicative generators seems to be known.

\medskip

This work was began during a staying at the Tata Institute of
Fundamental Research in January 2004. It is a pleasure to thank the
members of TIFR for their friendly hospitality. The financial support
from CEFIPRA is also gratefully acknowledged.

\section{Multiplicity-free subvarieties of $G$-varieties}

\subsection{}

First we introduce notation and conventions on algebraic varieties and
linear algebraic groups; as general references, we will use the books 
\cite{Ha77} and \cite{Ja03}. The ground field $k$ is algebraically
closed, of arbitrary characteristic. By a variety, we mean a separated
integral scheme of finite type over $k$; by a subvariety, we mean a
closed subvariety. A variety provided with an action of a linear
algebraic group $G$ is called a $G$-variety.

A quasi-coherent sheaf $\cF$ on a $G$-variety $X$ is $G$-linearized if
the symmetric algebra $S(\cF)$ is provided with an action of $G$,
compatible with the grading of this algebra and with the $G$-action on
$X$. The $G$-linearized quasi-coherent sheaves and their morphisms
form an abelian category. Those which are locally free are precisely 
the sheaves of sections of $G$-equivariant vector bundles on $X$. If 
$X$ is isomorphic to a homogeneous space $G/K$, where $K$ is a closed
subgroup of $G$, then each $K$-module $M$ defines a $G$-linearized
quasi-coherent sheaf $\cL_{G/K}(M)$ on $X$, with fiber $M$ at the base
point $K$. This yields an equivalence between the categories of
$K$-modules and of $G$-linearized quasi-coherent sheaves on $G/K$; the
finite-dimensional $K$-modules correspond to the $G$-linearized
coherent sheaves.

Now let $G$ be a connected reductive algebraic group. Let $B$ be a
Borel subgroup of $G$, and $T$ a maximal torus of $B$; let $B^+$ be
the opposite Borel subgroup, i.e., $B^+\cap B = T$. Let $U$, $U^+$ be
the unipotent radical of $B$, resp.~$B^+$. Denote by $W=N_G(T)/T$ the
Weyl group and by $R$ the set of roots of $(G,T)$. Let $R^+\subseteq
R$ be the subset of positive roots, i.e., of roots of $(B^+,T)$, and
let $S\subseteq R^+$ be the corresponding subset of simple roots. 
(Note that the roots of $(B,T)$ are negative). The weight lattice is
the character group of $T$, denoted by $\Lambda$; we identify this
group with the character group of $B$. The subset of dominant weights
associated with $R^+$ is denoted by $\Lambda^+$.

Each root $\alpha$ defines a reflection $s_\alpha \in W$; we may
identify the set $S$ with the subset of $W$ consisting of simple
reflections. These generate the group $W$; the length $\ell(w)$ of any
$w\in W$ is the minimal length of a word expressing $w$ in terms of
the simple reflections. There exists a unique element of maximal
length in $W$, denoted by $w_0$. For any simple reflection $s$, we
have a minimal parabolic subgroup $P_s: = B \cup BsB$ containing
$B$. The homogeneous space $P_s/B$ is isomorphic to the projective
line $\bP^1$. 

The flag variety of $G$ is the homogeneous space $G/B$, the disjoint
union of the Bruhat cells $C_w: = BwB/B$ ($w\in W$). Each $C_w$ is a
locally closed subvariety of dimension $\ell(w)$. Its closure in
$G/B$ is the Schubert variety $S_w = \overline{BwB}/B$, the union of
the cells $C_x$, where $x\in W$ and $x\le w$ for the Bruhat order.

Next let $X$ be a $G$-variety and let $Y$ be a
$B$-stable subvariety. The group $B$ acts on $G\times Y$ by 
$b (g,y) = (gb^{-1},by)$, and the quotient is a variety denoted by
$G\times^B Y$. The action of $G$ on itself by left multiplication
yields an action on $G\times^B Y$. The map 
$$
f=f_Y: G\times^B Y\to G/B, \quad (g,y)B\mapsto gB
$$ 
is a $G$-equivariant locally trivial fibration with fiber $Y$. 
On the other hand, the map 
$$
\pi=\pi_Y:G\times^B Y\to X, \quad (g,y)B\mapsto g y
$$
is well-defined, and the product map 
$$
f\times \pi: G\times^B Y \to G/B\times X
$$
is a $G$-equivariant closed immersion. Thus, $\pi$ is a proper
$G$-equivariant morphism, so that its image $GY$ is closed in $X$.
We will henceforth assume that $GY=X$. 

The fibers of $\pi$ may be identified with closed subschemes of $G/B$
via $f$. Specifically, the (scheme-theoretic) fiber $F_x$ at any point
$x\in X$ is isomorphic to the subscheme 
$\{gB~\vert~ g^{-1}x\in Y\}$ of $G/B$ 
(the image in $G/B$ of the preimage of $Y$ under the morphism 
$G\to X$, $g\mapsto g^{-1}x$). This subscheme is closed and stable
under the isotropy subgroup $G_x$.

More generally, for any $w\in W$, we have $B$-equivariant morphisms
$$
f_{Y,w}: \overline{BwB}\times^B Y \to \overline{BwB}/B = S_w
$$
(a locally trivial fibration with fiber $Y$), and
$$
\pi_{Y,w}: \overline{BwB}\times^B Y \to \overline{BwY}
$$
(a proper, surjective morphism). The fiber of $\pi_{Y,w}$ at any
$x\in \overline{BwY}$ may be identified with the scheme-theoretic
intersection $F_x\cap S_w$ in $G/B$. 

In particular, for any mimimal parabolic subgroup $P=P_s$, 
we have $P$-equivariant morphisms 
$f_{Y,s}: P\times^B Y \to P/B \cong \bP^1$ and
$\pi_{Y,s}: P\times^B Y \to PY$.

\medskip

\noindent
{\it Examples.} (1) If $Y=X$, then $G\times^B Y$ may be identified
with $G/B\times X$ so that $f$, $\pi$ are the projections.

\smallskip

\noindent
(2) If $X=G/B$, then $Y$ is a Schubert variety $S_w$. Further, $\pi$
is a $G$-equivariant locally trivial fibration with fiber $S_{w^{-1}}$ 
at the base point $B$.

\smallskip

\noindent
(3) As in the Introduction, assume that $Y$ contains a dense $B$-orbit 
$Y_0$ which is closed in the $G$-orbit $GY_0=:X_0$. Assume, in
addition, that the orbit maps $B\to Y_0$, $b\mapsto bx_0$ and 
$G\to X_0$, $g\mapsto gx_0$ are both separable for some $x_0\in Y_0$,
and the isotropy group $K:=G_{x_0}$ is connected. Put 
$B_K := B\cap K$, then $X_0\cong G/K$, $Y_0\cong B/B_K$, and
$\pi_0$ may be identified with the natural map $G/B_K \to G/K$. Since
$\pi_0$ is proper, it follows that $B_K$ is a Borel subgroup of
$K$. In other words, the fiber of $\pi_0$ at $x_0$ is the flag variety
of $K$.

\subsection{}

Let $X$ be a $G$-variety, and $Y$ a $B$-stable subvariety such that
$GY=X$. Following \cite{Br03}, we say that $Y$ is 
{\it multiplicity-free} if it satisfies the following two conditions:

\smallskip

\noindent
(i) Either $Y=X$, or $Y$ contains no $G$-orbit.

\smallskip

\noindent
(ii) For any minimal parabolic subgroup $P=P_s$ such that 
$P Y\neq Y$, the morphism $\pi_{Y,s}: P\times^B Y \to P Y$
is birational, and $P Y$ is multiplicity-free.

\medskip

Since $P Y$ is then a $B$-stable subvariety of $X$ of dimension
$\dim(Y)+1$, this definition makes sense by decreasing induction 
on the dimension of $Y$, starting with $Y=X$.

The multiplicity-free subvarieties of a nonsingular $G$-variety enjoy
remarkable properties: they are normal and Cohen-Macaulay by
\cite[Thm.1.2]{Br03}. If char$(k)=0$, then they have rational
singularities by \cite[Rem.3.3]{Br03}, generalizing
\cite[Thm.5]{Br01}. The Schubert varieties are examples of
multiplicity-free subvarieties.

We will obtain further homological properties of multiplicity-free
subvarieties; to state them, we introduce the following definition.

\medskip

We say that a proper morphism of varieties $\varphi: Z' \to Z$ is
$\cO$-{\it acyclic}, if the induced map 
$\cO_Z \to \varphi_*\cO_{Z'}$ is an isomorphism, and moreover 
$R^i\varphi_*\cO_{Z'} =0$ for all $i\ge 1$.

\medskip

(Such a morphism is also called {\it trivial} or {\it rational}, but
the name of $\cO$-acyclic is more specific). Note that any
$\cO$-acyclic morphism has connected fibers. Further, given two proper
morphisms of varieties $\psi:Z''\to Z'$ and  $\varphi: Z' \to Z$, the
composition $\phi \circ \psi$ is $\cO$-acyclic if both $\varphi$ and
$\psi$ are $\cO$-acyclic (as follows easily from the Leray spectral
sequence). Also, if both $\psi$ and $\phi \circ \psi$ are
$\cO$-acyclic, then $\varphi$ is $\cO$-acyclic. 

\begin{lemma}\label{tri}
Let $Y$ be a multiplicity-free subvariety of a complete nonsingular
$G$-variety $X$ and let $w\in W$. Then:

\noindent
(i) $\overline{BwY}$ is multiplicity-free. Further, $\pi_{Y,w}$ is
$\cO$-acyclic.

\noindent
(ii) If char$(k)=0$ and $B$ has an open orbit in $Y$, then
$H^i(\overline{BwY},\cO_{\overline{BwY}})=0$ for all $i\ge 1$. 
\end{lemma}

\begin{proof}
(i) We argue by induction on $\ell(w)$, the case where $\ell(w)=0$
being obvious. If $\ell(w)=1$, then $w$ is a simple reflection, say
$s$; also, put $P:=P_s$. Then $\overline{BsY} = P Y$ 
is multiplicity-free by definition. Thus, the morphism $\pi_{Y,s}$ is 
either birational or a fibration with fiber $\bP^1$. Since its image
$P Y$ is normal, it follows that 
$(\pi_{Y,s})_*\cO_{P\times^B Y}=\cO_{PY}$. Further, 
$R^i(\pi_{Y,s})_*\cO_{P\times^B Y}=0$ for all $i\ge 1$ 
by \cite[Lem.2.2]{Br03}. So $\pi_{Y,s}$ is $\cO$-acyclic.

If $\ell(w)\ge 2$, then we may write $w=vs$, where $v\in W$, $s\in S$,
and $\ell(w)=\ell(v)+1$. We then put $P:=P_s$ and $Z:=\overline{BwY}$,
so that $\overline{BwB}=\overline{BvB}P$, and $Z=\overline{BvP Y}$.
Since $PY$ is multiplicity-free, then so is $Z$ by the induction
assumption. Now consider the commutative diagram 
$$
\CD
\overline{BvB}\times^B P\times^B Y @>{\varphi}>> 
\overline{BwB}\times^B Y\\
@V{\psi}VV @V{\pi_{Y,w}}VV\\
\overline{BvB}\times^B PY 
@>{\pi_{P Y,v}}>> \overline{BwY},\\
\endCD
$$
where $\varphi((g,p,y)(B\times B))=(gp,y)B$ and
$\psi((g,p,y)(B\times B))= (g,py)B$.
Since $P Y$ is multiplicity-free, the morphism $\pi_{P Y,v}$
is $\cO$-acyclic by the induction assumption. On the other hand, we
have a cartesian diagram 
$$
\CD
\overline{BvB}\times^B P\times^B Y @>{\varphi}>> 
\overline{BwB}\times^B Y\\
@V{f_{Y,s}^v}VV  @V{f_{Y,w}}VV\\
\overline{BvB}\times^B P/B 
@>{f_{P/B,v}}>> S_w,\\
\endCD
$$
where $f_{Y,s}^v((g,p,y)(B\times B))=(g,p)B$. Further, $f_{Y,w}$ is a
locally trivial fibration, and the morphism $f_{P/B,v}$ is
$\cO$-acyclic (e.g., by \cite[Lem.2.2]{Br03} again). Thus, $\varphi$
is $\cO$-acyclic as well. Likewise, $\psi$ is $\cO$-acyclic. Since
$\pi_{Y,w} \circ \varphi = \pi_{PY,v} \circ \psi$, 
it follows that $\pi_{Y,w}$ is $\cO$-acyclic.

(ii) The variety $Y$ is rational, as it contains an open orbit of
the connected solvable linear algebraic group $B$. Since $P Y$ is 
birationally isomorphic to $Y$ or to $\bP^1\times Y$ via $\pi_{Y,s}$, 
it is rational as well. By induction on $\ell(w)$ again, it follows 
that $Z=\overline{BwY}$ is rational. Further, since $Z$ is 
multiplicity-free, it has rational singularities, i.e., there exists 
a $\cO$-acyclic desingularization $\varphi: Z'\to Z$. 
Thus, the Leray spectral sequence yields isomorphisms 
$H^i(Z,\cO_Z) \cong H^i(Z',\cO_{Z'})$. Now the nonsingular variety
$Z'$ is complete and rational, since so is $Z$. Thus,
$H^i(Z',\cO_{Z'})=0$ for all $i\ge 1$, by \cite[Lem.1]{Se59}. 
\end{proof}

\subsection{}

We now come to our main geometric result, which generalizes
\cite[Thm.6]{Br01} (specifically, the latter concerns spherical
varieties in characteristic $0$). The ingredients of the proof are
already in \cite{Br01} and \cite{Br03}; we will give details for
completeness.

\begin{theorem}\label{fib}
Let $Y$ be a multiplicity-free subvariety of a complete nonsingular
$G$-variety $X$. Then the morphism $\pi:G\times^B Y\to X$ is flat and
its fibers are connected, reduced and Cohen-Macaulay. In particular,
the general fibers of $\pi$ are varieties, and the fiber at any 
$B$-fixed point is a reduced union of Schubert varieties. 
\end{theorem}

\begin{proof} We divide it into five steps.

\smallskip

\noindent
{\it Step 1.} $Y$ intersects properly any $G$-orbit closure in $X$.

We check this by decreasing induction on the dimension of $Y$,
the case where $Y=X$ being trivial. 
Let $\cO$ be a $G$-orbit in $X$ and let $C$ be an irreducible 
component of $Y\cap \overline{\cO}$. Then $C\ne\overline{\cO}$, since
$Y$ contains no $G$-orbit. Thus, we may choose a minimal parabolic
subgroup $P \supset B$ such that $P C\neq C$. Then $P Y\neq Y$,
as $(PY)\cap\overline{\cO} = P (Y\cap\overline{\cO})$. By the
induction assumption, 
$$
\dim(P C) \le \dim(PY) + \dim(\cO) - \dim(X).
$$
But $\dim(P C) = \dim(C) + 1$ and $\dim(P Y) = \dim(Y) + 1$, so that 
$$
\dim(C) \le \dim(Y) + \dim(\cO)-\dim(X).
$$
On the other hand, the opposite inequality holds since $X$ is
nonsingular.

\smallskip

\noindent
{\it Step 2}. $\pi$ is equidimensional.

Indeed, its fiber $F_x$ at any $x\in X$ satisfies
$$\displaylines{
\dim(F_x) = -\dim(B) + \dim(\{g\in G~\vert~ g x\in Y\})
\hfill\cr
= -\dim(B) + \dim (Y\cap G x)+\dim(G_x)
\hfill\cr
= -\dim(B) + \dim(Y)+\dim(G x) -\dim(X) +\dim(G_x)
\hfill\cr
= -\dim(B) + \dim(Y) + \dim(G) -\dim(X)
= \dim(G\times^B Y) - \dim(X),
\hfill}$$
where the third equality follows from Step 3.

\smallskip

\noindent
{\it Step 3.} $\pi$ is flat and its fibers are Cohen-Macaulay.

Indeed, $G\times^B Y$ is Cohen-Macaulay, since $Y$ is. Further, $X$ is
nonsingular, and $\pi$ is equidimensional. Now the assertion follows
from \cite[Exer.III.10.9]{Ha77}.

\smallskip

\noindent
{\it Step 4.} The fiber $F_0$ at any $B$-fixed point is a reduced
union of Schubert varieties.

Indeed, since $F_0$ is $B$-stable, its irreducible components are
Schubert varieties. Thus, it suffices to show that $F_0$ is reduced. 
Since it is Cohen-Macaulay, it suffices in turn to show that it is
generically reduced, i.e., its multiplicities along irreducible
components are all $1$. This is equivalent to showing that the
intersection number $\int_{G/B} F_0 \cdot S_w$ is either $0$ or $1$,
for any $w\in W$ such that $\ell(w) + \dim(F_0) = \dim(G/B)$, i.e.,
$\ell(w) + \dim(Y) = \dim(X)$.

Now if $\overline{BwY}=X$, then the assumption on $w$ implies that the
morphism $\pi_{Y,w}$ is generically finite. Together with Lemma
\ref{tri} (i), it follows that $\pi_{Y,w}$ is birational. Thus, 
the scheme-theoretic intersection $F_x\cap S_w = \pi_{Y,w}^{-1}(x)$ is
a unique (reduced) point, for general $x\in X$. Since the number
$\int_{G/B} F_x \cdot S_w$ is independent of the point $x\in X$ by
flatness of $\pi$, this number must be $1$. On the other hand, if
$\overline{BwY} \neq X$ then $F_x\cap S_w$ is empty for general 
$x\in X$, and hence $\int_{G/B} F_x \cdot S_w=0$.

\smallskip

\noindent
{\it Step 5.} The fibers of $\pi$ are connected and reduced.

Indeed, the connectedness follows from Lemma \ref{tri}, and the
reducedness from Steps 3 and 4.
\end{proof}

\subsection{}

We still consider a complete nonsingular $G$-variety $X$ and a
multiplicity-free subvariety $Y$. Let $Z$ be another $B$-stable
subvariety of $X$ such that $GZ=X$. We will study the intersection 
$gY\cap Z$, where $g\in G$ is arbitrary. Clearly, this intersection
only depends on the double coset $BgB$. Thus, we may assume that $g\in
N_G(T)$ and we write $wY\cap Z$ for $gY\cap Z$, where $w=gT\in W$.

\begin{proposition}\label{int}
The set $\{w\in W ~\vert~ wY \text{ meets } Z \}$
has a unique maximal element $w(Y,Z)=w$ for the Bruhat order.
The scheme-theoretic intersection $wY\cap Z$ is a variety. Further,
$vY\cap Z$ is nonempty and connected for any $v\in W$, $v\leq w$. 
\end{proposition}

\begin{proof}
Let $V$ be the preimage of $Z$ under the morphism $\pi$. Consider the
diagram
$$
\CD
G/B @<{\varphi}<< V @>{\psi}>> Z\\
@V{\id}VV @V{\iota}VV @V{i}VV\\
G/B @<{f}<< G\times^B Y @>{\pi}>> X\\
\endCD
$$
where $i$ and $\iota$ are the inclusion maps, the square on the right 
is cartesian, and the square on the left is commutative.
Then $V$ is a closed $B$-stable subscheme of $G\times^B Y$, and 
$\varphi$, $\psi$ are $B$-equivariant. The fiber of $\varphi$ 
at any coset $gB$ may be identified with the intersection $gY\cap Z$.
On the other hand, by Theorem \ref{fib}, $\pi$, and hence $\psi$, is 
flat with connected and reduced fibers. Since $GZ =X$ and the general
fibers of $\pi$ are varieties, then so are the general fibers of
$\psi$. Together with the flatness of $\psi$, this implies that $V$ is
a variety.

Thus, $\varphi(V)$ is a $B$-stable subvariety of $G/B$,
i.e., a Schubert variety $S_w$. Since the Bruhat cell
$C_w$ is isomorphic to $U\cap w U w^{-1}$ via the
orbit map $g\mapsto gB$, and since $\varphi$ is equivariant,
the map 
$$
(U\cap w U w^{-1}) \times \varphi^{-1}(wB) \to V, \quad
(g,z)\mapsto g z
$$
is an open immersion with image $\varphi^{-1}(C_w)$. Thus, all the
fibers of $\varphi$ over $C_w$ are varieties. Since $S_w$ is normal,
Zariski's Main Theorem implies that $\varphi_*\cO_V = \cO_{S_w}$, so
that all the fibers of $\varphi$ are nonempty and connected.
\end{proof}

We say that $w(Y,Z)$ is the {\it generic position of $Y$ with respect
to $Z$}.

\medskip

For example, if $X=G/B$ then $Y$, $Z$ are Schubert varieties
$S_y$, $S_z$, where $y,z\in W$. Then $wY$ meets $Z$ if and only if 
there exists $y'\in W$ such that $y'\le y$ and $wy' \le z$. By
Proposition \ref{int}, the set of all such $w$ is an interval for the
Bruhat order. This may be seen directly: let $z':=wy'$, then 
$w=z'y'^{-1}$ is in the closure of $BzBy^{-1}B$. But this closure 
contains a unique open double class, $Bw(Y,Z)B$.

\section{A construction of $G$-linearized sheaves}

\subsection{}

Let $X$ be a $G$-variety, and $Y$ a $B$-subvariety such that $GY=X$. 
With the notation of 1.1, any $B$-module $M$ yields a locally free
$G$-linearized sheaf $\cL(M):=\cL_{G/B}(M)$ on $G/B$. This defines, 
in turn, a locally free $G$-linearized sheaf $f^*\cL(M)$ on 
$G\times^B Y$, and hence $G$-linearized sheaves 
$$
R^i\Phi_Y(M):= R^i \pi_*(f^*\cL(M)) \quad (i\ge 0)
$$
that we denote by $R^i\Phi(M)$ if this yields no confusion.
These sheaves are quasi-coherent; if, in addition, $M$ is
finite-dimensional, then they are coherent, as $\pi$ is proper. Also,
note that the sheaf $\Phi(M) := R^0\Phi(M) = \pi_*(f^*\cL(M))$ is
torsion-free, as $\pi$ is surjective.  

This defines additive functors $R^i\Phi$ from the category
of $B$-modules (resp.~finite dimensional $B$-modules) to the category
of quasi-coherent (resp. coherent) $G$-linearized sheaves on $X$. 
Since the fibers of $\pi$ are closed subschemes of $G/B$, the functors
$R^i\Phi$ vanish for all $i>\dim(G/B)$.
Also, note the isomorphism of $G$-linearized sheaves
$$
R^i\Phi(M\otimes N) \cong R^i\Phi(M) \otimes N, 
$$
where $M$ is any $B$-module, and $N$ is any $G$-module.

Clearly, the functor $\Phi$ is left exact. In fact, its $i$-th right
derived functor is $R^i\Phi$ as follows from the next result.  

\begin{lemma}
Let $I$ be an injective $B$-module. Then the sheaf $\Phi(I)$ is
generated by its global sections. Further, $R^i\Phi(I)=0$ for 
all $i\ge 1$.
\end{lemma}

\begin{proof}
We follow the argument as in the proof of \cite[Prop.1.1.4(i)]{Po89},
where the case of $X=G/B$ is treated. Recall that any injective
$B$-module is a direct factor of a direct sum of copies of the algebra
of regular functions $k[B]$, where $B$ acts by left multiplication. 
Thus, we may assume that $I=k[B]$. Now consider the diagram 
$$
\CD
G\times Y @>{\rho}>> G\times^B Y @>{\pi}>> X,\\
@V{\varphi}VV @V{f}VV\\
G @>{r}>> G/B
\endCD
$$
where $r$ denotes the quotient map, and the square is cartesian. 
Then the morphisms $r$ and $\rho$ are affine. Further, 
$\psi: = \pi \circ \rho$ is also affine, since it factors as the 
closed immersion 
$$
j : G\times Y \to G\times X, \quad (g,y)\mapsto (g,g y)
$$
followed by the projection 
$$
p: G\times X\to X.
$$ 
Clearly, $\cL(k[B]) = r_*\cO_G$, and hence
$$
R^i\Phi(k[B]) = R^i\pi_*(f^* r_* \cO_G)
= R^i\pi_*(\rho_*\varphi^* \cO_G) 
= R^i\pi_*(\rho_*\cO_{G\times Y}) 
= R^i\psi_*(\cO_{G\times Y}),
$$
where the second equality follows from the fact that cohomology
commutes with flat base extension \cite[Prop.III.9.3]{Ha77},
and the fourth equality from the exactness of $\rho_*$. 
Since $\psi$ is affine, this yields the vanishing of
$R^i\Phi(k[B])$ for all $i\ge 1$. Further, the sheaf
$$
\Phi(k[B]) = \psi_*(\cO_{G\times Y})
= p_* j_*(\cO_{G\times Y})
$$
is a quotient of $p_*\cO_{G\times X}=k[G]\otimes\cO_X$. In other
words, $\Phi(k[B])$ is generated by its subspace $k[G]$ of global
sections.
\end{proof}

Likewise, we define functors 
$$
R^i\Psi_Y := (R^i f_*) \circ \pi^* \quad (i\ge 0)
$$
from the category of quasi-coherent $G$-linearized sheaves on $X$, to
the corresponding category on $G/B$. We will identify the latter
category to the category of $B$-modules via $M\mapsto \cL(M)$. We put
$\Psi_Y := R^0\Psi = f_* \circ \pi^*$, and abbreviate $R^i\Psi_Y$, 
resp.~$\Psi_Y$ by $R^i\Psi$, resp.~$\Psi$ if this yields no confusion. 
We may now state

\begin{lemma}\label{adj} 
(i) $R^i\Psi = H^i(Y,-)$ for any $i\ge 0$.

\noindent
(ii) For any locally free $G$-linearized sheaf $\cF$ on $X$, and for
any $B$-module $M$, we have an isomorphism of $G$-modules
$$
\Gamma(X,\Phi(M)\otimes\cF) \cong 
\Gamma(G/B,\cL(M\otimes \Psi(\cF))).
$$
In particular, $\Gamma(X,\Phi(M)) \cong \Gamma(G/B,\cL(M))$ if $X$ is 
complete.
\end{lemma}

\begin{proof}
(i) Consider a $G$-linearized sheaf $\cG$ on $G\times^B Y$. Since
$f:G\times^B Y\to G/B$ is a locally trivial fibration with fiber $Y$,
then 
$$
R^if_*(\cG)= \cL(H^i(Y,\iota^*\cG)),
$$ 
where $\iota:Y\to G\times^B Y$ denotes the inclusion. Now take 
$\cG = \pi^*\cF$, where $\cF$ is a $G$-linearized sheaf on $X$. Then
$\iota^*\cG$ is the restriction of $\cF$ to $Y$, since $\pi\circ\iota$
is the inclusion of $Y$ into $X$.

\noindent
(ii) follows readily from the projection formula.
\end{proof}

Note that the $G$-module $\Gamma(G/B,\cL(M))$ is the induced module
$\Ind_B^G(M)$. Further, $H^i(G/B,\cL(M))=R^i\Ind_B^G(M)$, where
$R^i\Ind_B^G$ denotes the $i$-th right derived functor of induction.

Next, for any $\lambda \in \Lambda$, we denote by $k_{\lambda}$
the one-dimensional $B$-module with weight $\lambda$,
and we put $\cL(\lambda):=\cL(k_\lambda)$. Then each $\cL(\lambda)$ is
a $G$-linearized invertible sheaf on $G/B$. Further, $\cL(\lambda)$ is
generated by its global sections if and only if 
$\lambda \in \Lambda^+$. Also put 
$$
\Phi(\lambda):=\Phi(k_{\lambda}), \quad 
R^i\Phi(\lambda):= R^i\Phi(k_{\lambda}).
$$
The isomorphisms 
$\cL(\lambda)\otimes \cL(\mu) \cong \cL(\lambda+\mu)$
for any $\lambda,\mu\in\Lambda$, yield isomorphisms
$f^*\cL(\lambda)\otimes f^*\cL(\mu) \cong f^*\cL(\lambda+\mu)$
and hence, morphisms of sheaves of $\cO_X$-modules
$$
\Phi(\lambda)\otimes \Phi(\mu) \to \Phi(\lambda+\mu).
$$ 
These morphisms define a $\Lambda$-graded algebra structure
on the sheaf $\bigoplus_{\lambda\in \Lambda} \Phi(\lambda)$.

By Lemma \ref{adj}, the $G$-module $\Gamma(X,\Phi(\lambda))$ is
isomorphic to the {\it dual Weyl module} 
$\Gamma(G/B,\cL(\lambda))=\Ind_B^G(\lambda)$,
if $X$ is complete. This yields an isomorphism of $\Lambda$-graded
rings
$$
\bigoplus_{\lambda\in \Lambda} \Gamma(X,\Phi(\lambda))
\cong \bigoplus_{\lambda\in \Lambda} \Ind_B^G(\lambda)
\cong \Gamma(G/U,\cO_{G/U}).
$$

We now describe the functors $R^i\Phi$, $R^i\Psi$ for the examples 
already considered in 1.1. 

\medskip

\noindent
{\it Examples.} (1) If $Y=X$, then 
$R^i\Phi_X(M) = H^i(G/B,\cL(M))\otimes \cO_X$
for any $B$-module $M$. In other words, 
$R^i\Phi_X = (R^i\Ind_B^G)\otimes \cO_X$.  

\smallskip

\noindent
(2) If $X=G/B$ and $Y=S_w$, then 
$$
R^i\Phi_{S_w}(M)= \cL(H^i(S_{w^{-1}},\cL(M))).
$$
In other words, $R^i\Phi_{S_w}$ (regarded as an endofunctor of the
category of $B$-modules) is the {\it Joseph functor} $H^i_{w^{-1}}$
introduced in \cite{Po89} (see also \cite{vdK93}). On the other hand,
$R^i\Psi_{S_w}$ is the Joseph functor $H^i_w$ by Lemma \ref{adj}
(i). The $B$-module $\Psi(\lambda) = \Gamma(S_w,\cL(\lambda))$,
where $\lambda \in \Lambda^+$ and $w\in W$, is the 
{\it dual Joseph module} denoted by $H_w(\lambda)$; further, 
$H^i_w(\lambda)=0$ for any $i\ge 1$, by \cite[Prop.1.4.2]{Po89}. In
particular, $H_{w_0}(\lambda)$ is the dual Weyl module
$\Ind_B^G(\lambda)$, and $R^i\Ind_B^G(\lambda)=0$ for any $i\ge 1$,
again for $\lambda \in \Lambda^+$. 

\smallskip

\noindent
(3) Let $X_0$ be an $G$-stable open subset of $X$ and put 
$Y_0:=Y\cap X_0$, then the restriction to $X_0$ of any sheaf
$R^i\Phi_Y(M)$ is $R^i\Phi_{Y_0}(M)$.  

Now assume that $\pi_0:B\times^B Y_0 \to X_0$ is the fibration
$G/B_K \to G/K$, where $K$ is a connected subgroup of $G$
with Borel subgroup $B_K := B\cap K$. Then the $G$-linearized sheaf
$R^i\Phi_{Y_0}(M)$ is associated with the $K$-module
$H^i(K/B_K,\cL(M))$. Thus, for any $\lambda \in \Lambda^+$, the
sheaves $R^i\Phi_{Y_0}(\lambda)$ vanish for $i\ge 1$, since
$\cL(\lambda)$ restricts to a globally generated invertible sheaf on
the flag variety $K/B_K$. Further, $\Phi_{Y_0}(\lambda)$ is
associated to $\Gamma(K/B_K,\cL(\lambda))$, a dual Weyl module for 
the quotient of $K$ by its radical. Since any endomorphism of a dual
Weyl module is a scalar, it follows that any endomorphism of the
$G$-linearized sheaf $\Phi_Y(\lambda)$ is a scalar as well. In
particular, $\Phi_Y(\lambda)$ is indecomposable as a $G$-linearized
sheaf. 

\subsection{}

We keep the notation and assumptions of 2.1 and we will study the 
maps in (equivariant) $K$-theory induced by the functors $R^i\Phi$.

Let $K^G(X)$ be the Grothendieck group of the abelian category 
of $G$-linearized coherent sheaves on $X$. Then $K^G(X)$ is a module
over the Grothendieck ring of finite-dimensional $G$-modules, i.e., 
over the representation ring $R(G)$. If, in addition, $X$ is 
nonsingular, then $K^G(X)$ is isomorphic to the Grothendieck group 
of $G$-linearized locally free sheaves on $X$, and hence is 
an $R(G)$-algebra.

For example, the Grothendieck ring $K^G(G/B)$ is isomorphic to 
the representation ring $R(B)\cong R(T)$. In turn, 
$R(T) \cong \bZ[\Lambda]$ (the group ring of $\Lambda$) and
this isomorphism identifiess $R(G)$ to $\bZ[\Lambda]^W$ 
(the subring of $W$-invariants).

For any $G$-linearized coherent sheaf $\cF$ on $X$, let 
$[\cF]$ be its class in $K^G(X)$. Then assigning to each 
finite-dimensional $B$-module $M$ the class 
$\sum_i (-1)^i \, [R^i\Phi(M)]$ yields a map 
$$
\Phi_!^G: R(T) \to K^G(X)
$$
which is $R(G)$-linear, but generally not a ring homomorphism, 
as shown by the examples below.

Likewise, we have the Grothendieck group $K(X)$ of coherent sheaves 
on $X$, and an additive map
$$
\Phi_!: K(G/B) \to K(X), \quad 
[\cF] \mapsto \pi_!(f^![\cF]) = \sum_i (-1)^i \, [R^i\pi_*(f^*\cF)].
$$
Clearly, $\Phi_!$ is compatible with $\Phi_!^G$ via the forgetful 
maps $R(T) \to K(G/B)$ (the characteristic map) and $K^G(X) \to K(X)$. 
Recall from \cite{Dem74} that the abelian group $K(G/B)$ is freely 
generated by the Schubert classes $[\cO_{S_w}]$, $w\in W$.
Further, the rational Grothendieck ring $K(G/B)_{\bQ}$ is isomorphic 
(via the Chern character) to the rational Chow ring $A^*(G/B)_{\bQ}$,
and the latter ring is generated by Chern classes of line 
bundles. It follows easily that the characteristic map is surjective 
over the rationals (it is surjective if $G$ is semisimple and
simply connected, see \cite{Ma76}). Thus, the vector space 
$K(G/B)_{\bQ}$ is spanned by the classes $[\cL(\lambda)]$, 
$\lambda\in\Lambda^+$.

As a consequence, the subgroup $\Phi_!(R(T))$ is generated 
by the $\pi_![\cO_{\overline{BwB}\times^B Y}]$, $w \in W$; 
this subgroup is generated over $\bQ$ by the $\Phi_!(\lambda)$, 
$\lambda\in\Lambda^+$.

\medskip

\noindent
{\it Examples.} (1) If $X$ is a point, then we obtain a surjective 
map
$$
\chi^G: R(T) \to R(G), \quad 
[M] \mapsto \sum_i (-1)^i \, [R^i\Ind_B^G(M)],
$$
given by Weyl's character formula. Thus, $\chi^G$
(regarded as a self-map of $\bZ[\Lambda]$) is the Demazure operator 
associated with $w_0$.

More generally, if $Y=X$, then $\Phi^G_!(z)= \chi^G(z) \, [\cO_X]$ 
for any $z\in R(T)$. Thus, the image of $\Phi_!^G$ is the 
$R(G)$-submodule of $K^G(X)$ generated by $[\cO_X]$.

\smallskip

\noindent
(2) If $X=G/B$ and $Y=S_w$, then 
$\Phi^G_!: \bZ[\Lambda] \to \bZ[\Lambda]$
is the Demazure operator associated with $w^{-1}$.

\smallskip

\noindent
(3) Assume that $X$ contains an open orbit $X_0 = G/K$, 
where $K$ is a connected subgroup of $G$ with Borel 
subgroup $B_K :=B\cap K$, and that $Y$ is the closure in $X$ of 
$Y_0 = BK/K \cong B/B_K$. Then we may assume that $T_K:=T\cap K$
is a maximal torus of $K$. The composition of $\Phi^G_!$ with 
the restriction map $K^G(X) \to K^G(X_0) \cong R(K)$ is the map 
$$
R(T) \to R(K), \quad [M]\mapsto \sum_i (-1)^i \, R^i\Ind_{B_K}^K(M),
$$
the composition of the restriction $R(T) \to R(T\cap K)$
with $\chi^K: R(T\cap K) \to R(K)$. As a consequence,
this composition is surjective.

\subsection{}

We now apply the results of 1.2 and 1.3 to study the sheaves
$R^i\Phi(\lambda)$ under the assumption of multiplicity-freeness. 

\begin{theorem}\label{loc}
Let $X$ be a complete nonsingular $G$-variety, $Y$ a multiplicity-free
subvariety, and $\lambda \in \Lambda^+$. Then:

\noindent
(i) The $G$-linearized sheaf $\Phi(\lambda)$ is locally free and 
generated by its space of global sections, $\Ind_B^G(\lambda)$. 
The fiber of $\Phi(\lambda)$ at any $x\in X$ equals
$\Gamma(F_x,\cL(\lambda))$ (where $F_x$ denotes the fiber of $\pi$ at
$x$, regarded as a closed subscheme of $G/B$). Further,
$R^i\Phi(\lambda) = 0$ for any $i\ge 1$.

\noindent
(ii) The morphism
$\Phi(\lambda) \otimes \Phi(\mu) \to \Phi(\lambda + \mu)$ 
is surjective for any $\mu \in \Lambda^+$.

\noindent
(iii) The image of the map $\Phi_!: K(G/B) \to K(X)$ is the subgroup
generated by the $[\cO_{\overline{BwY}}]$, $w\in W$; this 
subgroup is generated over $\bQ$ by the $[\Phi(\lambda)]$, 
$\lambda \in \Lambda^+$.

\noindent
(iv) If char$(k)=0$ and $B$ has an open orbit in $Y$,
then $H^i(X,\Phi(\lambda))=0$ for any $i\ge 1$.
\end{theorem} 

\begin{proof}
(i) Recall from Theorem \ref{fib} that the fiber $F_0$ at any
$B$-fixed point is a reduced union of Schubert varieties. By
\cite{Ra87}, it follows that the restriction map 
$\Gamma(G/B,\cL(\lambda)) \to \Gamma(F_0,\cL(\lambda))$
is surjective, and $H^i(F_0,\cL(\lambda))=0$ for all $i\ge 1$. Since
$\pi$ is flat and $B$-equivariant, and any $B$-orbit closure in $X$
contains a fixed point, these properties extend to arbitrary fibers by
semicontinuity \cite[Thm.III.12.8]{Ha77}. Our assertions now
follow by cohomology and base change \cite[Thm.III.12.11]{Ha77}.

(ii) follows as above from the surjectivity of the product map
$$
\Gamma(F_0,\cL(\lambda)) \otimes \Gamma(F_0,\cL(\mu))
\to \Gamma(F_0,\cL(\lambda+\mu)),
$$
established in \cite{Ra87}.

(iii) follows from the discussion in 2.2 together with (i) and 
Lemma \ref{tri} (i).

(iv) Since $R^j\pi_*(f^*\cL(\lambda)) = 0$ for any $j\ge 1$, 
the Leray spectral sequence for $\pi$ yields isomorphisms for all $i$ 
$$
H^i(X,\Phi(\lambda)) \cong H^i(G\times^B Y,f^*\cL(\lambda)).
$$
To compute the latter, we use the Leray spectral sequence for
$f$. For any $j\ge 0$, the $G$-linearized sheaf on $G/B$
$$
R^jf_*(f^*\cL(\lambda)) = 
\cL(\lambda)\otimes R^jf_*\cO_{G\times^B Y}
$$
corresponds to the $B$-module $k_{\lambda}\otimes H^j(Y,\cO_Y)$.
By Lemma \ref{tri}, this module is zero for all $j\ge 1$, if
char$(k)=0$ and $Y$ contains an open $B$-orbit. Thus,
$$
H^i(G\times^B Y,f^*\cL(\lambda))\cong H^i(G/B,\cL(\lambda)).
$$
But the latter vanishes for all $i\ge 1$, since 
$\lambda \in \Lambda^+$; this completes the proof.
\end{proof}

\subsection{} 

We still consider a complete nonsingular $G$-variety $X$ and a
multiplicity-free subvariety $Y$. We will obtain several relations
between the functors $R^i\Phi_Y$, $R^i\Psi_Y$ and the Joseph functors
$H^i_w$. We begin with the following  

\begin{proposition}
There are spectral sequences
$$
R^i\Phi_Y(H^j_{w^{-1}}(M)) \Rightarrow R^{i+j}\Phi_{\overline{BwY}}(M)
$$
and 
$$
H^i_w(R^j\Psi_Y(\cF)) \Rightarrow R^{i+j}\Psi_{\overline{BwY}}(\cF),
$$
where $w\in W$, $M$ is a $B$-module, and $\cF$ is a locally free
$G$-linearized sheaf on $X$.
\end{proposition}

\begin{proof}
Put $Z:=\overline{BwY}$. The morphism 
$\pi_{Y,w}:\overline{BwB}\times^B Y\to Z$ induces a morphism
$$
\rho:G\times^B \overline{BwB}\times^B Y \to G\times^B Z, \quad
(g,h,y)(B\times B)\mapsto (g,hy)B
$$
which is $\cO$-acyclic, since $\pi_{Y,w}$ is (by Lemma \ref{tri}). 
Together with the projection formula, this implies isomorphisms
$$
R^i\Phi_Z(M) = R^i\pi_{Z*}(f_Z^*\cL(M)) = 
R^i(\pi_Z\circ \rho)_*(\rho^*f_Z^*\cL(M)).
$$
On the other hand, $\pi_Z\circ \rho = \pi_Y\circ\psi$, where $\psi$
denotes the morphism
$$
G\times^B \overline{BwB}\times^B Y \to G\times^B Y, \quad
(g,h,y)(B\times B)\mapsto (gh,y)B.
$$
This yields a spectral sequence
$$
R^i\pi_{Y*}(R^j\psi_*(\rho^*f_Z^*\cL(M)))\Rightarrow
R^{i+j}(\pi_Z\circ \rho)_*(\rho^*f_Z^*\cL(M)))=R^i\Phi_Z(M). 
$$
To obtain the first spectral sequence, it suffices to show that
$R^j\psi_*(\rho^*f_Z^*\cL(M))= f_Y^*\cL(H^j_{w^{-1}}(M))$.
For this, consider the commutative diagram
$$
\CD
G/B @<{\varphi}<< G\times^B \overline{BwB}\times^B Y @>{\psi}>>
G\times^B Y \\
@V{\id}VV @V{p}VV @V{f_Y}VV \\
G/B @<{f_w}<< G\times^B S_w @>{\pi_w}>> G/B,\\
\endCD
$$
where $p((g,h,y)(B\times B))=(g,hB)B$ and 
$\varphi((g,h,y)(B\times B))=gB$, so that $\varphi=f_Z\circ \rho$. 
Note that the square on the right is is a cartesian square of locally
trivial fibrations. This yields
$$\displaylines{
R^j\psi_*(\rho^*f_Z^*\cL(M)) = R^j\psi_*(\varphi^*\cL(M))
=R^j\psi_*(p^*f_w^*\cL(M)) 
\cr = f_Y^*R^j\pi_{w*}(f_w^*\cL(M)) = f_Y^*\cL(H^j_{w^{-1}}(M))
}$$
as desired, where the third equality follows again from 
\cite[Prop.III.9.3]{Ha77}.

To obtain the second exact sequence, note that
$$
R^i\Psi_Z(\cF) = H^i(Z,\cF) = 
H^i(\overline{BwB}\times^B Y, \pi_{Y,w}^* \cF)
$$
by Lemma \ref{adj} (i), Lemma \ref{tri} (i) and the projection
formula. Thus, the morphism $f_{Y,w}:\overline{BwB}\times^B Y \to Z$
yields a Leray spectral sequence 
$$
H^i(S_w,R^jf_{Y,w*}(\pi_{Y,w}^* \cF)) \Rightarrow R^{i+j}\Psi_Z(\cF).
$$
To complete the proof, it suffices to show that the $B$-linearized sheaf
$R^jf_{Y,w*}(\pi_{Y,w}^*\cF)$ is isomorphic to the restriction of
$R^jf_*(\pi^*\cF)$ to $S_w$. For this, consider the cartesian square
$$
\CD
\overline{BwB}\times^B Y @>{f_{Y,w}}>> S_w \\
@V{\iota}VV @V{i}VV \\
G\times^B Y @>{f}>> G/B,\\
\endCD
$$
where the vertical arrows are the inclusions. By
\cite[Rem.III.9.3.1]{Ha77}, there is a natural map
$i^*R^jf_*(\pi^*\cF) \to R^jf_{Y,w*}(\iota^*\pi^*\cF)$.
Further, $\iota^*\pi^*\cF =\pi_{Y,w}^*\cF$, as 
$\pi_{Y,w} = \pi \circ \iota$. This yields a morphism of
$B$-linearized sheaves on $S_w$, 
$$
\phi:i^*R^jf_*(\pi^*\cF) \to R^jf_{Y,w*}(\pi_{Y,w}^*\cF).
$$
Further, the sheaf $R^jf_*(\pi^*\cF) = \cL(H^j(Y,\cF))$ has fiber
$H^j(gY,\cF)$ at any $gB \in G/B$; and, since $f_{Y,w}$ is a locally
trivial fibration with fibers being the translates $gY$,
$gB \in S_w$, the fiber of $R^jf_{Y,w*}(\pi_{Y,w}^*\cF)$ at $gB$ also
equals $H^j(gY,\cF)$. Finally, the induced map $\phi_{gB}$ on fibers
is the identity, so that $\phi$ is an isomorphism.
\end{proof}

Taking $M = k_{\lambda}$ where $\lambda\in\Lambda^+$, and using 
the vanishing of $H^i_{w^{-1}}(\lambda)$ for all $i\ge 1$ together
with Theorem \ref{loc}, this implies readily the following result.

\begin{corollary}\label{jos}
For any $\lambda\in\Lambda^+$ and $w\in W$, we have
$\Phi_Y(H_w(\lambda))= \Phi_{\overline{Bw^{-1}Y}}(\lambda)$
and $R^i\Phi_Y(H_w(\lambda))=0$ for all $i\ge 1$.
\end{corollary}

Next recall from \cite{Po89}, \cite{vdK93} that a finite-dimensional
$B$-module $M$ is said to admit an {\it excellent filtration}, if there
exists a sequence $0=M_0\subset M_1\subset\cdots\subset M_n = M$ of 
$B$-submodules, such that each subquotient $M_i/M_{i-1}$ is a dual
Joseph module. Also, a finite-dimensional $G$-module $N$ is said to
admit a {\it good filtration}, if there exists a sequence
$0=N_0\subset N_1\subset\cdots\subset N_n = N$ of $G$-submodules,
such that each subquotient $N_i/N_{i-1}$ is a dual Weyl module.
Now Corollary \ref{jos} and Theorem \ref{loc} imply readily

\begin{corollary}\label{exc}
Let $M$ be a finite-dimensional $B$-module admitting an excellent
filtration. Then the sheaf $\Phi_Y(M)$ is locally free and generated
by its global sections, and $R^i\Phi_Y(M) = 0$ for all $i\ge 1$.
Further, the $G$-module $\Gamma(X,\Phi_Y(M))$ admits a good filtration. 

If char$(k)=0$ and $B$ has an open orbit in $Y$, then 
$H^i(X,\Phi_Y(M))=0$ for all $i\ge 1$.
\end{corollary}

Finally, we obtain one additional relation between the functors
$\Phi_Y$, $\Psi_Z$ and $H_w$. 

\begin{proposition}\label{com}
Let $Z$ be a $B$-stable subvariety of $X$ such that $G Z=X$
and let $w$ be the generic position of $Y$ with respect to $Z$
(as defined in 1.4). Then $\Psi_Z(\Phi_Y(\lambda))= H_w(\lambda)$
for any $\lambda \in \Lambda^+$.
\end{proposition}

\begin{proof}
By Lemma \ref{adj}, 
$\Psi_Z(\Phi_Y(\lambda)) = \Gamma(Z,\Phi_Y(\lambda))$.
We now use the notation of the proof of Proposition \ref{int}.
By \cite[Rem.III.9.3.1]{Ha77}, there is a natural map
$$
\varepsilon : i^*\Phi_Y(\lambda) = 
i^*\pi_*f^*\cL(\lambda) \to \psi_*\iota^*f^*\cL(\lambda).
= \psi_*\varphi^*\cL(\lambda).
$$
Further, arguing as in the proof of Theorem \ref{loc}, one checks that
the right-hand side is a locally free sheaf on $Z$ and that
$\varepsilon$ restricts to an isomorphism on all the fibers. 
Thus, $\varepsilon$ is an isomorphism. Taking global 
sections yields 
$$
\Gamma(Z,\Phi_Y(\lambda)) = \Gamma(V,\varphi^*\cL(\lambda))
= \Gamma(G/B,\cL(\lambda)\otimes \varphi_*\cO_V)
= \Gamma(S_w,\cL(\lambda))= H_w(\lambda),
$$
since $\varphi_*\cO_V = \cO_{S_w}$ by (the proof of) Proposition
\ref{int}. 
\end{proof}

\section{Symmetric spaces of minimal rank}

\subsection{} 

We assume from now on that char$(k)\neq 2$. Let $\theta$ be an
automorphism of order $2$ of the group $G$. Denote by $K=G^{\theta}$
the subgroup of $\theta$-fixed points, then the homogeneous space
$G/K$ is a {\it symmetric space}.

It is known that the group $K$ is reductive; further, its identity
component $K^0$ is nontrivial unless $G$ is a torus and $\theta$ is
the inversion. We refer to \cite[Sec.2]{Sp85} for these facts and for
further results on symmetric spaces, mentioned in this subsection and in
the next one.

A $\theta$-stable subtorus $S$ of $G$ is said to be 
$\theta$-{\it fixed} (resp. $\theta$-{\it split}) 
if $\theta(x)=x$ (resp. $\theta(x)=x^{-1}$)
for all $x\in S$. On the other hand, a parabolic subgroup $P$ 
of $G$ is said to be $\theta$-{\it split}, if the parabolic
subgroup $\theta(P)$ is opposite to $P$, i.e., $P \cap \theta(P)$ 
is a Levi subgroup of both. 

The maximal $\theta$-fixed subtori of $G$ are just
the maximal tori of $K$. Thus, any two such tori are conjugate by an
element of $K^0$; their common dimension is the rank of $K$, denoted
by $\rk(K)$. Further, the centralizer in $G$ of any $\theta$-fixed
subtorus is a $\theta$-stable maximal torus of $G$. 

The maximal $\theta$-split subtori are also conjugate in $K^0$; 
their common dimension is called the {\it rank} of the symmetric space
$G/K$ and is denoted by $\rk(G/K)$. Further, any maximal torus of $G$
containing a maximal $\theta$-split subtorus $A$ is $\theta$-stable, and
any two such maximal tori are conjugate in $K^0$. The centralizer 
$C_G(A)$ is a Levi subgroup of a minimal $\theta$-split parabolic 
subgroup $P$, and the derived subgroup of $C_G(A)$ is contained in $K^0$.
Any two minimal $\theta$-split parabolic subgroups are conjugate in $K^0$.

Let $T$ be a $\theta$-stable maximal subtorus of $G$ and put
$T^{\theta} := \{ x\in T ~\vert~ \theta(x) = x \}$,
$T^{-\theta} := \{ x\in T ~\vert~ \theta(x) = x^{-1} \}$.
Then $T = T^{\theta} T^{-\theta}$ and $T^{\theta} \cap T^{-\theta}$
is finite, so that $\dim(T) = \dim(T^{\theta}) + \dim(T^{-\theta})$.
Since $\dim(T) = \rk(G)$, $\dim(T^{\theta}) \le \rk(K)$ and 
$\dim(T^{-\theta})\le \rk(G/K)$, this yields
$$
\rk(G/K) \ge \rk(G) - \rk(K).
$$

We say that the symmetric space $G/K$ {\it is of minimal rank}, if
$\rk(G/K) = \rk(G) - \rk(K)$. The preceding discussion yields several
equivalent formulations of this definition. 

\begin{lemma}\label{cond}
For the symmetric space $G/K$, the following conditions are equivalent:

\noindent
(i) $G/K$ is of minimal rank.

\noindent
(ii) Any $\theta$-stable maximal torus of $G$ contains a maximal
$\theta$-fixed subtorus.

\noindent
(iii) Any $\theta$-stable maximal torus of $G$ contains a maximal
$\theta$-split subtorus.

\noindent
(iv) Any two $\theta$-stable maximal tori of $G$ are conjugate in
$K^0$. 
\end{lemma}

\subsection{}

We still consider an involution $\theta$ of $G$ and the corresponding
symmetric space $G/K$. We may then choose a $\theta$-stable Borel
subgroup $B$ and a $\theta$-stable maximal torus $T$ of $B$; such
a pair $(T,B)$ is called a {\it standard pair}. Then $U$, $B^+$ and
$U^+$ are $\theta$-stable as well. Further, $\theta$ acts on $W$ and
on $R$; it stabilizes $R^+$, $R^-$ and $S$.

Denote by $\fg$, $\fb$, $\fu$, $\ft$, $\ldots$ the Lie algebras of
$G$, $B$, $U$, $T$, $\ldots$. Then $\theta$ acts on $\fg$ by an
automorphism of order $2$, still denoted by $\theta$. Further, 
$\fg^{\theta}=\fk$ and the decomposition 
$\fg = \fu \oplus \ft \oplus \fu^+$ is $\theta$-stable, whence 
$\fk = \fu^{\theta} \oplus \ft^{\theta} \oplus (\fu^+)^{\theta}$. 
It follows that the connected component $T^{\theta,0}$ is a maximal 
$\theta$-fixed torus, with centralizer $T$ in $G$, and that 
$B^{\theta,0}$, $(B^+)^{\theta,0}$ are opposite Borel subgroups 
of $K$. As a consequence, the orbit $B/(B\cap K)$ is closed in 
$G/K$.

We may now formulate several special properties of symmetric spaces of
minimal rank.

\begin{lemma}\label{comb}
Let $G/K$ be a symmetric space of minimal rank. Then:

\noindent
(i) We have the decomposition 
$$
\fk = \ft^{\theta} \oplus 
\bigoplus_{\alpha\in R^{\theta}} \fg_\alpha \oplus
\bigoplus_{\alpha\in R \setminus R^{\theta}}
k(x_\alpha + \theta(x_\alpha)),
$$
where for any root $\alpha$, we denote by $x_\alpha$ a generator of
the root subspace $\fg_\alpha\subset \fg$. Further, $R^{\theta}$ is 
the root system of the Levi subgroup containing $T$ of a minimal 
$\theta$-split parabolic subgroup.

\noindent
(ii) The roots of $(K^0,T^{\theta,0})$ are exactly the restrictions to 
$T^{\theta,0}$ of the roots of $(G,T)$. 

\noindent
(iii) The Weyl group $W_{K^0} := W(K^0,T^{\theta,0})$ may be identified 
with $W^{\theta}$; it contains $w_0$.

\noindent
(iv) Any standard pair is conjugate to $(T,B)$ by an element of
$K^0$.

If $G$ is semisimple and adjoint, then the groups 
$T_K:=T\cap K = T^{\theta}$ and $B_K:=B\cap K = B^{\theta}$ are
connected. Further, $K$ is connected, semisimple and adjoint. Finally,
any dominant weight of $K$ (with respect to its maximal torus $T_K$
and its Borel subgroup $B_K$) extends to a dominant weight of $G$
(w.r.t. $T$, $B$).
\end{lemma}

\begin{proof}

(i) Note that $R^{\theta}$ is the root system of $C_G(T^{-\theta,0})$,
where $T^{-\theta,0}$ is a maximal $\theta$-split torus.
Thus, if $\alpha\in R$ satisfies $\theta(\alpha)=\alpha$, then 
$\theta$ fixes pointwise $\fg_\alpha$. On the other hand, if 
$\theta(\alpha)\ne\alpha$, then $\fg_\alpha\oplus\fg_{\theta(\alpha)}$ 
is $\theta$-stable, and its $\theta$-fixed subspace is spanned by
$x_\alpha + \theta(x_\alpha)$. This implies our assertions. 

(ii) is a direct consequence of (i).

(iii) Let $x\in K^0$ normalize $T^{\theta,0}$, then $x$ normalizes 
its centralizer $T$ in $G$. This identifies $W_{K^0}$ with a 
subgroup of $W^{\theta}$. 

Conversely, let $w\in W^{\theta}$, then $wBw^{-1}$ is a 
$\theta$-stable Borel subgroup of $G$ containing $T$. Thus,
$wBw^{-1}$ contains a Borel subgroup of $K^0$ containing 
$T^{\theta,0}$. Replacing $w$ with $wx$ for some 
$x \in W(K^0,T^{\theta,0})$, we may assume that $wBw^{-1}$ 
contains $B^{\theta,0}$. But (i) implies that $B$ is
the unique Borel subgroup containing $B^{\theta,0}$.
It follows that $wBw^{-1}=B$, i.e., $w=1$.

So $W_{K^0}$ equals $W^{\theta}$. But the latter contains
$w_0$: indeed, $\theta w_0 \theta\in W$ sends $R^+$ to $R^-$, 
so that $\theta w_0 \theta = w_0$.

(iv) Let $(T_1,B_1)$ be a standard pair. By the results in 3.1, we
may assume that $T_1=T$ after conjugation by some element of
$K^0$. Then $B_1=wBw^{-1}$ for a unique $w\in W$. Further, $w$ has
a representative in $K^0$ by the argument of (iii).

Assume that $G$ is semisimple and adjoint, i.e., the simple
roots form a basis of the weight lattice $\Lambda$. Let $x\in T$,
then $x\in K$ if and only if $\alpha(\theta(x)) = \alpha(x)$
for all $\alpha\in S$. This is equivalent to 
$(\theta(\alpha)\alpha^{-1})(x)=1$ for all $\alpha\in S$. In other
words, $T_K$ is the intersection of the kernels of the
characters $\theta(\alpha)-\alpha$ of $T$, where $\alpha$ runs over
the nontrivial orbits of $\theta$ in $S$. Now these characters form
part of a basis of $\Lambda$, and hence $T_K$ is connected.
Since $U_K := U\cap K$ is connected by \cite[Prop.4.8]{Sp85}, 
it follows that $B_K = T_K U_K$ is connected as well. 

Consider the action of $K$ on the set of standard pairs by
conjugation, then the isotropy group of $(T,B)$ is $T_K$.
Together with (iv), it follows that $K= K^0 T_K$. Thus, $K$ is 
connected. So, by (ii), the positive roots of $(K,T_K)$ are
exactly the restrictions of the positive roots of $(G,T)$. It follows
that any dominant weight of $K$ extends to a dominant weight of $G$,
and that the intersection of the kernels of the roots of 
$(K,T_K)$ is trivial. Thus, $K$ is semisimple and adjoint.
\end{proof}

\subsection{}

Any symmetric space $G/K$ contains only finitely many $B$-orbits; 
they have been studied in detail in \cite{Sp85}, \cite{RiSp90}. 
By specializing these results to spaces of minimal rank, we obtain:

\begin{proposition}\label{borb}
Let $X_0 := G/K$ be a symmetric space of minimal rank and 
let $(T,B)$ be a standard pair. Then:

\noindent
(i) The assignement $w\in W\mapsto BwK/K$ induces a bijection from
$W/W_K$ to the set of $B$-orbits in $X_0$.

\noindent
(ii) Any $v \in W/W_K$ admits a representative 
$w\in W$ such that $\overline{BvK}/K=\overline{BwY_0}$.

\noindent
(iii) $Y_0$ is the unique closed $B$-orbit in $G/K$. Any
$B$-orbit closure is multiplicity-free.

\noindent
(iv) Any two $B$-orbit closures are in generic position $w_0$ 
(in the sense of 1.4).  
\end{proposition}

\begin{proof}
(i) Let $\cN:=\{x\in G ~\vert~ x\theta(x^{-1})\in N_G(T)\}$.
This is a closed subset of $G$, stable under the action of 
$T\times K$ via $(t,k)x = txk^{-1}$. By \cite[Thm.4.2]{Sp85}, 
any double coset $BxK$ in $G$ meets $\cN$ into a unique orbit 
of $T\times K$. On the other hand, $x^{-1}Tx$ is a $\theta$-stable
maximal torus for any $x\in\cN$. Thus, there exists $k\in K$
such that $x^{-1}Tx = k^{-1}Tk$, i.e., $\cN=N_G(T) K$. It follows
that 
$$
\cN/(T\times K) \cong N_G(T)/TN_K(T) \cong W/W_K.
$$

(ii) Let $w\in W$, then $x := w\theta(w^{-1})$ is an element 
of $W$ which only depends on the coset $wW_K$. Clearly, 
$\theta(x)=x^{-1}$, i.e., $x$ is a twisted involution in the 
sense of \cite[3.1]{Sp85}. Further, by Lemma \ref{comb} (i), 
any root in $R$ is either complex or compact imaginary 
relative to $x$, in the terminology of [loc.cit., 3.6]. 
So, by [loc.cit., Prop.3.3, Lem.3.7], we may write
$$
x = s_1 \cdots s_h \theta(s_h) \cdots \theta(s_1),
$$
where $s_1,\ldots,s_h$ are simple reflections, and $\ell(x) = 2h$. 
Then $w \in s_1 \ldots s_h W^{\theta}$. By Lemma \ref{comb} (iii),
we may assume that $w = s_1 \ldots s_h$; then this decomposition 
is reduced. So, by \cite[Thm.6.5]{Sp85}, 
$\overline{BwK}/K = \overline{BwB}Y_0 = \overline{BwY_0}$.

(iii) By (ii), any $B$-orbit closure contains $Y_0$. 
(Alternatively, by \cite[Cor.6.6]{Sp85}, the closed $B$-orbits 
are in bijection with the $K$-conjugacy classes of standard pairs. 
But there is a unique such class by Lemma \ref{comb} (iv).)
 
The assertion on multiplicity-freeness follows from the results
in \cite{RiSp90}. Specifically, since the roots with respect 
to any $B$-orbit are either complex or compact imaginary,
only cases I and II occur in the discussion of 
\cite[4.2, 4.3]{RiSp90}. It follows that the morphism 
$P\times^B Y\to PY$ is birational for any $B$-orbit closure $Y$
and for any minimal parabolic subgroup $P$ such that $PY \ne Y$.

(iv) It suffices to show that $Y_0\cong BK/K$ meets
$w_0Y_0\cong w_0BK/K$. But $w_0TK =TK$, since $w_0$ has a
representative in $N_K(T)$ by Lemma \ref{comb} (iii).
\end{proof}

\subsection{}

We assume from now on that $G$ is semisimple and adjoint. Then, 
by \cite{DePr83} in characteristic $0$ and \cite{DeSp99} in an 
arbitrary characteristic, any symmetric space $G/K$ admits a 
wonderful compactification $X$. This is a projective 
nonsingular $G$-variety containing an open orbit isomorphic 
to $G/K$. The complement of this open orbit is a union of 
nonsingular prime divisors with normal crossings, the
boundary divisors $X_1,\ldots,X_r$, where $r : = \rk(G/K)$. 
The $G$-orbit closures are exactly the partial intersections 
of these divisors. Finally, there is a unique closed orbit, 
the intersection of all the boundary divisors. This orbit is
isomorphic to $G/P$, where $P$ is a minimal $\theta$-split 
parabolic subgroup.

We now gather some more or less well-known facts on $B$-orbit 
closures in wonderful compactifications, which we will prove 
for completeness.

\begin{lemma}\label{wond}
Let $X$ be the wonderful compactification of a symmetric space 
$G/K$. Then:

\noindent
(i) The closure in $X$ of any non-open $B$-orbit in $G/K$ contains 
no $G$-orbit. 

\noindent
(ii) The abelian group $K(X)$ admits a basis indexed by $X^T$ 
(the set of $T$-fixed points in $X$), and this set is finite.
The classes $[\cO_Z]$, where $Z$ runs over the closures in $X$ 
of $B$-orbit in $G/K$, form part of such a basis.

\noindent
(iii) $G/K$ is of minimal rank if and only if the closed $G$-orbit
in $X$ contains $X^T$.
\end{lemma}

\begin{proof} 
(i) The open $B$-orbit being affine, its complement is of pure 
codimension $1$ in $G/K$. Thus, any non-open $B$-orbit in $G/K$ 
is contained in the closure of a $B$-orbit of codimension $1$. 
Let $Y$ be the closure in $X$ of this last orbit; then $Y$ does 
not contain the closed $G$-orbit, by \cite[p.~293]{DeSp99}.

(ii) Since $X$ contains only finitely many $G$-orbits, the set $X^T$ 
is finite. Then a suitable Bialynicki-Birula decomposition yields 
a paving of $X$ by locally closed $B$-stable cells indexed by $X^T$, 
each cell being isomorphic to some affine space; see \cite{DePr83}, 
\cite{DeSp85}. It follows that the classes $[\cO_Z]$ of the structure 
sheaves of cell closures form a basis of the group $K(X)$.
Further, by \cite[Thm.4.2]{DeSp85}, each cell intersects each $G$-orbit 
into a unique $B$-orbit. Thus, the closure in $X$ of any $B$-orbit in 
$G/K$ is the closure of some cell.

(iii) If $G/K$ is of minimal rank, then it contains no $T$-fixed 
points unless the group $G$ is trivial. Together with the 
description of isotropy groups of $G$-orbits in $X$ given
in \cite[Thm.5.2]{DePr83} (see also \cite{DeSp85}), it follows 
that only the closed orbit contains fixed points.

For the converse, let $Y_0$ be a closed $B$-orbit in $X_0=G/K$; then
$Y_0$ is contained in a unique cell $C$. Let $x$ be the unique
$T$-fixed point in $C$. By assumption, $x$ lies in the closed 
$G$-orbit, so that the boundary divisors $X_1,\ldots,X_r$ admit 
local equations at $x$, eigenvectors of $T$ with linearly independent 
weights $\chi_1,\ldots,\chi_r$. Therefore, $\chi_1,\ldots,\chi_r$ are 
the weights of $T$ in the normal space to $C\cap Gx = Bx$ in $C$ at $x$. 
Since $C$ is $T$-equivariantly isomorphic to its tangent space at $x$,
it follows that the subgroup of $T$ where 
$\chi_1 = \cdots = \chi_r = 1$ fixes points of $X_0 \cap C = Y_0$. 
Thus, $\dim(T) - r \ge \rk(K)$, that is, 
$\rk(G/K) \le \rk(G) - \rk(K)$. This forces the equality.
\end{proof}

Next we combine all the previous results to obtain our main
statement:

\begin{theorem}\label{ext}
Let $X$ be the wonderful compactification of a symmetric space 
$G/K$ of minimal rank. Choose a standard pair $(T,B)$ and put
$B_K:=B\cap K$, $T_K:=T\cap K$; let $Y$ be the closure in $X$ of the
closed $B$-orbit $B/B_K$ in $G/K$. Finally, let $\mu$ be a dominant
weight of $K$ (with respect to $T_K$, $B_K$) and choose a dominant 
weight $\lambda$ of $G$ (w.r.t. $T$, $B$) which extends $\mu$. Then:

\noindent
(i) The $G$-linearized sheaf $\Phi(\lambda)$ on $X$ is locally free 
and generated by its space of global sections, $\Ind_B^G(\lambda)$. 
The restriction of $\Phi(\lambda)$ to $G/K$ is the $G$-linearized 
sheaf associated with the $K$-module $\Ind_{B_K}^K(\mu)$.
Further, the restriction map 
$\Gamma(X,\Phi(\lambda)) \to \Gamma(Y,\Phi(\lambda))$
is an isomorphism.

\noindent
(ii) The subspace of $K(X)_{\bQ}$ spanned by the classes 
$[\Phi(\lambda)]$, $\lambda\in\Lambda^+$, has basis the classes 
$[\cO_Z]$, where $Z$ runs over the closures of $B$-orbits in $G/K$.

\noindent
(iii) If char$(k)=0$, then $H^i(X,\Phi(\lambda))=0$ for all $i\ge 1$.
\end{theorem}

\begin{proof}
(i) $Y$ is multiplicity-free by Proposition \ref{borb} (iv) and
Lemma \ref{wond} (i). So Theorem \ref{loc} applies to yield
the first assertion. The second assertion holds by Example 3 in
2.1. The third assertion is a consequence of Lemma \ref{adj} (i), 
Proposition \ref{com} and Proposition \ref{borb} (iv).

(ii) and (iii) follow from Theorem \ref{loc} together with Lemma 
\ref{wond} (ii). 
\end{proof}

This statement implies all the assertions of Theorem \ref{main},
except for the uniqueness of the extension $\Phi(\lambda)$. 
But this is a consequence of the following observation, a variant 
of \cite[Prop.1]{BoSe58}.

\begin{lemma}\label{unique}
Let $X_0$ be an open subset of a variety $X$. Let $\cF_0$ be 
a torsion-free, coherent sheaf on $X_0$, generated by a 
finite-dimensional subspace $V$ of $\Gamma(X_0,\cF_0)$.
Then there exists a unique torsion-free coherent sheaf $\cF$ 
on $X$ extending $\cF_0$ and generated by its subspace $V$ of 
global sections.
\end{lemma}

\begin{proof}
Let $\cG_0$ be the kernel of the evaluation map 
$V \otimes \cO_{X_0} \to \cF_0$ and let $\cG$ be the subsheaf of 
$V \otimes \cO_X$ consisting of those local sections $s$ such that 
$s\vert_{X_0}$ lies in $\cG_0$. Clearly, $\cG$ is a coherent sheaf 
on $X$ which extends $\cG_0$. This yields an exact sequence of 
sheaves on $X$:  
$$
0 \to \cG \to V \otimes \cO_X \to \cF \to 0,
$$
where $\cF$ is a coherent sheaf on $X$ which extends $\cF_0$.
Further, if $f$ is a nonzero local section of $\cO_X$, and $s$ is
a local section of $V \otimes \cO_X$ such that $f s \in \cG$, then
$s\in \cG$ since $\cF_0$ is torsion-free. Thus, $\cF$ is 
torsion-free as well.

If $\cF'$ is another sheaf satisfying the assertions of the 
lemma, then we have an exact sequence
$$
0 \to \cG' \to V \otimes \cO_X \to \cF' \to 0.
$$
Moreover, $\cG'$ is contained in $\cG$, since $\cF'$ extends
$\cF_0$. This yields an exact sequence
$$
0 \to \cT \to \cF' \to \cF \to 0,
$$
where $\cT$ is a coherent sheaf supported on $X\setminus X_0$. 
Since $\cF'$ is torsion-free, $\cT=0$, whence $\cF' = \cF$.
\end{proof}

\subsection{}

We conclude with some comments on the assertions of Theorem 
\ref{ext}.

\smallskip

\noindent
(i) The fibers of the equivariant vector bundle $\Phi(\lambda)$ 
are the quotients of the vector space $\Ind_B^G(\lambda)$ 
obtained by composing the (surjective) restriction map 
$$
r : \Ind_B^G(\lambda) \to \Ind_{B_K}^K(\mu)
$$
with elements of $G$ and taking limits in the Grassmannian. 
Note that the $K$-module $\Ind_B^G(\lambda)$ admits a good 
filtration with head $\Ind_{B_K}^K(\mu)$; see \cite[Sec.3]{Bru97} 
and \cite[Thm.29]{vdK01}.

\smallskip

\noindent
(ii) By Proposition \ref{borb}, the number of $B$-orbits in $G/K$ 
equals the index $[W:W_K] = [W:W^{\theta}]$. On the other hand, 
by Lemmas \ref{comb} (i) and \ref{wond}, the dimension of 
$K(X)_{\bQ}$ equals the index $[W:W(R^{\theta})]$. Further,
$W(R^{\theta})$ is a proper subgroup of $W^{\theta}$. So we see 
that the classes $[\Phi(\lambda)]$, $\lambda\in\Lambda^+$, generate
a proper subspace of $K(X)_{\bQ}$. But one may show that 
$K(X)_{\bQ}$ is a free module of rank $[W^{\theta}:W(R^{\theta}]$
over its subring generated by the classes 
$[\cO_{X_1}],\ldots,[\cO_{X_r}]$ of the boundary divisors, 
together with the classes 
$[\Phi(\lambda_1)],\ldots,[\Phi(\lambda_s)]$, where 
$\lambda_1,\ldots,\lambda_s$ are extensions of the 
$s=\rk(K)$ fundamental weights of $K$. Further, these 
$r+s =\rk(G)$ classes are algebraically independent. This will 
be developed elsewhere.

\smallskip

\noindent
(iii) The vanishing of the higher cohomology groups 
$H^i(X,\Phi(\lambda))$ in an arbitrary characteristic 
would follow from the vanishing of the groups $H^j(Y,\cO_Y)$, 
by the argument of Theorem \ref{loc}. 
In the case of the group $K\times K/\diag(K)$, the latter vanishing
holds by \cite[Cor.3]{BrPo00}. It would be interesting to know if 
this holds in our setting of symmetric spaces of minimal rank. 
More generally, given the closure $Z$ of a $B$-orbit in $G/K$, and 
a globally generated $G$-equivariant line bundle $\cL$ on $X$, 
it would be interesting to know if the $B$-module $H^0(Z,\cL)$ 
admits a Schubert filtration, and if $H^j(Z,\cL)=0$ whenever $j\ge 1$. 
This holds in the group case by \cite[Cor.3,Thm.7]{BrPo00}.

\end{document}